\documentclass[twoside]{amsart}
\usepackage{amssymb,graphicx,rotating,multirow,multicol}
\usepackage{amsmath,amsthm,amsfonts,amscd}
\usepackage[mathscr]{euscript}

\swapnumbers

{\catcode`\@=11 \gdef\mnote#1{\marginpar{\footnotesize
 \tolerance\@M\spaceskip2.6\p@ plus10\p@ minus.9\p@\rm#1}}}

\let\Bbb\mathbb

\def\Z{\Bbb Z}
\def\R{\Bbb R}

\let\+\sqcup
\def\fac#1#2{{#1/\lower0.5ex\hbox{\small{#2}}}}
\def\ztwo{\Z/\hspace*{-0.4ex}\lower0.5ex\hbox{\tiny2}}
\def\sfit{\sffamily\itshape}

\makeatletter
\newcommand{\addresseshere}{%
  \enddoc@text\let\enddoc@text\relax
}
\def\section{\@startsection{section}{1}%
  \normalparindent{.5\linespacing\@plus.7\linespacing}{-.5em}%
  {\normalfont\bfseries\Large}}
\makeatother

\title{Rokhlin's signature theorems}
\author{S.~Finashin, V.~Kharlamov, O.~Viro }
\begin{document}
\maketitle
\thispagestyle{empty}

\maketitle

\noindent{\bf The scope of Rokhlin's mathematics.}
Rokhlin's mathematical heritage includes a number of outstanding
contributions
to various fields:  measure and ergodic theories,
topology, and real algebraic geometry. 
He used to change the focus of his research interests in
average once in five years.
In each of the periods, he
managed to obtain fundamental results and then
 turned to a new field, being attracted by a challenging problem.

Rokhlin started in the measure and ergodic theories. The results of his PhD 
were presented in his paper ``On the fundamental concepts of measure theory''.
The second (doctoral) dissertation defended soon after his PhD was called 
``On the most important metric classes of dynamical systems''. The first
ergodic theory period finished in 1950. In the late fifties, 
Rokhlin came back to ergodic theory in his research related to 
entropy of dynamical systems.  In the first half of the fifties, the 
period of topology followed. 
\vspace{5pt}

\noindent{\bf A good challenge.}
When Rokhlin turned to topology, his  initial challenge was to find the 
third stable homotopy group of spheres, that is $\pi_{n+3}(S^n)$ with 
$n\ge5$. Prior to that Pontryagin found the groups 
$\pi_{n+2}(S^n)$ with $n\ge2$ and $\pi_{n+1}(S^n)$ with $n\ge3$ by 
his method, which bridges the homotopy theory and differential topology.
Technically, Pontryagin's calculation was a study of curves and surfaces 
embedded in Euclidean space with trivialized normal bundles.   
It relied on a comparatively simple topology of curves
and surfaces. 
A topology of 3- and 4-dimensional smooth manifolds that was
needed for a similar calculation of $\pi_{n+3}(S^n)$ had not yet
developed. \vspace{5pt}

\noindent{\bf Groundbreaking results.} In a striking tour de force, 
Rokhlin developed a technique, which allowed him to prove that
$\pi_{n+3}(S^n)=\Z/_{24}$ for $n\ge5$, and found himself in a new research
area, where he could make groundbreaking discoveries. Below we overview
the results published by Rokhlin in the notes of 1951-52 (see \cite{R1,R2,R3,R4}), where he  
calculated $\pi_{n+3}(S^n)$. 
The results grew out of calculations of cobordism groups $\Omega_3$ and $\Omega_4$. 
(Cobordisms and the cobordism groups $\Omega_n$ and
$\mathfrak N_n$ were introduced by Rokhlin in the same notes.)
\vspace{2pt} 

\noindent$\bullet$ {\sfit Any oriented closed smooth 3-manifold bounds a 
compact oriented 4-manifold.\/} (In other words,  $\Omega_3=0$.) 
\vspace{2pt} 

Rokhlin introduced the signature $\sigma(M)$ of an oriented closed
$4n$-manifold $M$. 
The signature $\sigma(M)$ is the difference between the numbers of positive
and negative coefficients in a diagonalization of the intersection form 
$H_{2n}(M;\R)\times H_{2n}(M;\R)\to\R$ of $M$. By discovering the
properties of signature stated below
Rokhlin turned it into one of the central invariants in topology of
manifolds.
\vspace{2pt} 
 
\noindent$\bullet$ 
{\sfit The signature of an oriented smooth closed 4n-manifold vanishes 
if the manifold bounds an oriented smooth compact $(4n+1)$-manifold.\/} 
(In other words, the signature defines 
a homomorphism $\Omega_{4n}\to\Z$.)
\vspace{2pt} 

\noindent$\bullet$ {\sfit An oriented smooth closed 4-manifold $M$ 
bounds a compact oriented 5-manifold if and only if $\sigma(M)=0$.
In general, an oriented smooth closed 4-manifold $M$ is cobordant 
to the disjoint sum 
of $\sigma(M)$ copies of the complex projective plane.}
(In particular,  $\Omega_4=\Z$.) 
\vspace{2pt} 

\noindent$\bullet$ {\sfit 
The signature of an oriented smooth closed 4-manifold $M$ equals one 
third of the Pontryagin number $p_1(M)[M]$.\/} 
\vspace{2pt} 
 
\noindent$\bullet$ {\sfit The signature of an oriented smooth closed 4-manifold $M$ 
with $w_2(M)=0$ is divisible by 16. }
\vspace{2pt} 

Most of these achievements are often attributed to other authors, who  
discovered these things later. However, the last item in this list stands
apart. It is commonly referred to as {\sfit the Rokhlin Theorem\/} and 
considered Rokhlin's  most famous result. Many times it played a
substantial r\^ole in the subsequent development of topology. Below we 
concentrate on a few of them that seem the most important. 
\vspace{5pt}

\noindent{\bf Reformulations of the Rokhlin Theorem.} To a non-specialist, the
Rokhlin theorem sounds a bit cryptic and technical. Let us take a closer look.
The theorem establishes a relation between 
two basic characteristics of an oriented smooth closed 4-manifold $M$: 
its second Stiefel-Whitney class $w_2(M)$ and the  
signature $\sigma(M)$. It claims that if $w_2(M)=0$, 
then $\sigma(M)\equiv0\mod 16$.

What does the assumption $w_2(M)=0$ mean? 
The Stiefel-Whitney classes 
$w_k(M)\in H^k(M;\ztwo)$
measure complexity 
of the tangent bundle $TM$. Orientability of $M$ means $w_1(M)=0$.
For an orientable smooth closed 4-manifold $M$,  $w_2(M)$ is the only
obstruction for $M$ being {\sfit almost parallelizable}, that is admitting 
4 tangent vector fields  
linear independent at all points but one. With this in mind we can 
reformulate the Rokhlin Theorem as follows:\vspace{2pt}

\noindent
$\bullet$ {\sfit The signature of an  almost parallelizable smooth closed
4-manifold is divisible by 16.}

 A smooth closed manifold $M$ admits a Spin-structure iff $w_1(M)=0$ and
$w_2(M)=0$. Thus we can reformulate the Rokhlin Theorem as follows: 
\vspace{2pt}

\noindent$\bullet$
{\sfit The signature of a smooth closed Spin 4-manifold is divisible by 16.}

For an orientable smooth closed 4-manifold $M$ the assumption $w_2(M)=0$ 
holds true iff the $\ztwo$
 intersection form 
$H_2(M;\ztwo)\times H_2(M;\ztwo)\to\ztwo$ is even, that is the
self-intersection number of any closed surface in $M$ is even. This
interpretation does not depend on the smooth structure or tangent bundle of
$M$. 
Thus we get one more reformulation:\\
\noindent$\bullet$
{\sfit If an oriented smooth closed 4-manifold $M$ has even
$\ztwo$ intersection form, then $\sigma(M)\equiv0\mod16$.}
\vspace{5pt}
 
\noindent{\bf Causes for divisibility: Topology versus Algebra.}
Divisibility of some integer by 16 looks quite
mysterious. Sixteen is a large non-prime number. As we will see, 
a half of divisibility of the signature by 16, that is divisibility 
by 8, can be extracted by purely algebraic arguments  
from simple topological facts, which are not specific for smooth 
4-manifolds. On the contrary, the last factor 2 in 16 manifests a crucial 
obstruction for a topological manifold to admit a smooth  structure. 

Due to Poincar\'e duality, the intersection form $H_2(M)\times H_2(M)\to\Z$
of a closed oriented 4-manifold $M$
is a symmetric integral unimodular form.  If $w_2(M)=0$, then this 
form is even, because its reduction modulo 2 is a part of 
the $\ztwo$ intersection form, which is even as $w_2(M)=0$.
The signature of an even symmetric integral unimodular form is 
divisible by 8.  This is a purely algebraic fact. 

The Rokhlin Theorem claims that $\sigma(M)$ is divisible not only by 8, 
as the algebra ensures, but by 16.
 So, it imposes a restriction on a symmetric even unimodular 
form realizable as the intersection form of an oriented smooth closed 
4-manifold $M$ with $w_2(M)=0$. Namely, the signature of a realizable 
form cannot 
be congruent to $8\mod16$. For example, the form $E_8$ (a famous even
unimodular integral form of rank and signature 8) is prohibited. 
\vspace{5pt}

\noindent{\bf Obstructions to Diff or PL on a 4-manifold.}
In Rokhlin's theorem, the assumption that $M$ is {\sfit smooth\/} is 
necessary: there exists an oriented closed
topological simply-connected 4-manifold with any unimodular intersection 
form, as it was proven by Freedman [1982]. 
If the form is even and its signature is not divisible by 16 
(say, if this is $E_8$) then by the Rokhlin Theorem the manifold 
does not admit a smooth structure. In lower dimensions this does not 
happen: any $n$-manifold  with $n<4$ admits a smooth structure.

In dimension four, smoothability is equivalent to existence of a piecewise 
linear structure (PL-structure). Thus the Rokhlin theorem provides an
obstruction to existence of a PL-structure. \vspace{5pt}

\noindent{\bf The key to PL in high dimensions.}
In the eighties Donaldson
discovered other numerous obstructions to existence of PL or smooth
structure on 4-manifolds. However,  they are less
robust. The Rokhlin Theorem is formulated in terms invariant under
cobordisms and gives rise to {\sfit high-dimensional\/} results, while
the other obstructions do not. 

For example, if
$M$ is a simply-connected closed 4-manifold with signature non-divisible by
16 (and does not admit a PL-structure by the Rokhlin Theorem), 
then $M\times (S^1)^n$ does not admit a PL-structure for any $n>0$.  This
can be deduced from the Rokhlin Theorem.

Complete obstructions to existence or equivalence of PL-structures on 
manifolds of dimension $\ge5$ (discovered by Kirby and Siebenmann in
the seventies) rely on the Rokhlin Theorem.       
\vspace{5pt}

\noindent{\bf More about the Rokhlin Theorem.}
Due to the space restrictions, our story is very incomplete.  
We said nothing about numerous
generalizations of the Rokhlin Theorem, their applications, proofs,
misconceptions, etc.
You can find all of this on 
http://mathcenter.spb.ru/nikaan/book/index.html . 
See also a book \cite{Book} by L. Guillou and A. Marin.

\end{document}